\documentclass[a4paper]{article}
\usepackage{amssymb}
\usepackage{amsmath}
\usepackage{amsfonts}

\newtheorem{theorem}{Theorem}[section]

\newtheorem{remark}[theorem]{Remark}

\newtheorem{conjecture}[theorem]{Conjecture}

\newtheorem{assumption}[theorem]{Assumption}

\newenvironment{Proof}{\removelastskip\par\medskip
\noindent{\em Proof.}
\rm}{\penalty-20\null\hfill$\square$\par\medbreak}

\newcommand{\Q}{\mathbb{Q}}

\newcommand{\F}{\mathbb{F}}

\begin{document}
\title{Modularity of rigid Calabi-Yau threefolds over $\Q$}
\author{Luis Dieulefait and Jayanta Manoharmayum}

\date{}
\maketitle

\begin{abstract}
We prove modularity for a huge class of rigid Calabi-Yau
threefolds over $\Q$. In particular we prove that every rigid
Calabi-Yau threefold with good reduction at $3$ and $7$ is
modular.
\end{abstract}

\section{Introduction}

Let $X$ be a rigid Calabi-Yau threefold defined over $\Q .$ Recall
that (see \cite{noriko}, \cite{noriko-update})  $X$ is a
smooth projective threefold over $\Q,$ and satisfies
\begin{enumerate}
\item $H^1(X,\mathcal{O}_X)=H^2(X,\mathcal{O}_X)=0,$
\item $\bigwedge^3\Omega_X^1\cong \mathcal{O}_X,$ and
\item $h^{2,1}(X)=0.$
\end{enumerate}
The first two of these conditions define a Calabi-Yau threefold,
while the third signifies rigidity. In a certain sense, rigid
Calabi-Yau varieties are the natural generalization of elliptic
curves to higher dimensions.

We approach the issue of modularity for rigid Calabi-Yau
threefolds by using Galois representations. For other
formulations, see \cite{noriko}, \cite{noriko-update}. Here is how
it works: Given a prime $\ell,$ the action of $G_\Q$ on the dual
of the \'{e}tale cohomology group $H^3_{\textrm{\'{e}t}}(X_{\bar\Q},
\Q_\ell)$ defines
  a continuous, odd, two dimensional representation
$$\rho_{X,\ell}:G_\Q \longrightarrow GL_2(\Q_\ell).$$
The representation $\rho_{X,\ell}$ is unramified at a prime not
equal to $\ell$ if  $X$ has good reduction there; its determinant
is equal to $\epsilon_\ell^3$ (where $\epsilon_\ell$ is the
$\ell$-adic cyclotomic character). If $X$ has good reduction at
$\ell$ and $\ell>3$, it follows from the comparison theorem of
Fontaine-Messing (cf. \cite{FM}) that the representation $\rho_{X,
\ell}$ is crystalline (locally at $\ell$). Moreover, from the
definition of rigid Calabi-Yau it also follows that this local
$\ell$-adic representation has Hodge-Tate weights $\{ 0 , 3 \}$.

If $X$ has good reduction at $p$, we set
$$t_3(p)=\textrm{tr}\rho_{X,\ell}(\textrm{Frob}_p),\quad \ell\neq
p.$$ Of course, the definition is independent of which $\ell$ we
choose for calculating the trace. Moreover, the quantity $t_3(p)$
is relatively easy to calculate: one only needs to `count points'
on the mod $p$ reduction of $X$ (see Remark 2.9 of
\cite{noriko-update} for details).

The Riemann hypothesis, proved by Deligne as part of the Weil
conjectures, states that the roots of the characteristic
polynomial of $\rho_{X, \ell}(\textrm{Frob}_p)$, $p \neq \ell$ a
prime of good reduction,  both have  the same complex absolute
value. If $\rho_{X, \ell}$ was reducible, its constituents should
be $\ell$-adic characters of $G_\Q$, and so we must have
$$\rho_{X,\ell}\sim \begin{pmatrix}\mu \epsilon_\ell^a & *\\
0&\mu^{-1} \epsilon_\ell^b\end{pmatrix} $$ where $\mu$ is a finite
order character and $a,b \in \mathbb{Z}$ with $a+b = 3.$   Since
this clearly  contradicts the Riemann hypothesis as $a \neq b,$ we
can therefore conclude that $\rho_{X,\ell}$ is irreducible.

\begin{conjecture}
Let $X$ be a rigid Calabi-Yau threefold over $\Q.$ There is then a
weight $4$ newform $f$ on $\Gamma_0(N)$ such that
$$\rho_{X,\ell}\sim\rho_{f,\ell}$$
for some (and hence every) prime $\ell.$
\end{conjecture}
We haven't specified what the level should be in the statement of
the conjecture, but it is clear that $N$ should only be divisible
by primes where $X$ has bad reduction. Specifying the level at
primes of bad reduction is more problematic: in fact, results of
de Jong (see in particular Deligne's corollary in \cite{Ber},
Prop.~6.3.2) imply that when $\ell$ ranges through all primes of
good reduction of $X$,
 there exists a  bound for the (prime to $\ell$ part of the)
 conductor of $ \rho_{X,\ell} $, and this minimal upper bound will turn out to be the
  level of the associated modular form. However, this value is not effectively
  computable.
   As for evidence in
support of the conjecture, it is known to be true for all rigid
Calabi-Yau threefolds  that we know of. (However, we must add that
`what we know' contains only a  handful of examples!) In the form
stated above, the conjecture forms a part of (or rather follows
from) the Fontaine-Mazur conjecture on modularity of Galois
representations.

In this article, we prove the
modularity conjecture for all rigid Calabi-Yau threefolds verifying some local
 conditions. More precisely, we shall prove the
following:
\begin{theorem}
Let $X$ be a rigid Calabi-Yau threefold over $\Q.$ Suppose $X$
satisfies one of the following two conditions:
\begin{enumerate}
\item $X$ has good reduction at $3$ and $7$; or,
\item $X$ has good reduction at $5$  and some prime $p\equiv \pm
2\pmod{5}$ with $t_3(p)$ not divisible by $5.$
\end{enumerate}
Then $X$ is modular.
\end{theorem}

\begin{remark}
 This result applies to all but one (class) of the known examples of rigid Calabi-Yau
threefolds.
\end{remark}

\noindent{\bf Acknowledgements:} We would like to thank Frazer
Jarvis and the referee for various helpful suggestions.

\section{Modularity of $\ell$-adic representations}

More generally, let us consider irreducible continuous
representations
$$\rho:G_\Q \longrightarrow GL_2(\Q_\ell^{ac})$$
which are odd and unramified outside finitely many primes. We also
assume that $\ell\geq 5$, and that the determinant of $\rho$ is
$\epsilon_\ell^3.$ By choosing a stable lattice, we can reduce
$\rho$ modulo $\ell$ and obtain a unique semi-simple continuous
representation
$$\overline\rho:G_\Q \longrightarrow GL_2(\F_\ell^{ac}).$$
Following Wiles, one then deduces the modularity of $\rho$ from
that of $\overline\rho.$

\subsection{}
We shall now describe some results under which such a deduction is
possible. We fix $\rho$ and $\overline\rho$ as above (in the
beginning of the main section). Throughout this subsection, we
shall also suppose that
\begin{assumption}
$\rho$ is crystalline at $\ell$ with Hodge-Tate weights $\{0 , 3\}$.
\end{assumption}
 The first result we describe is due to Skinner and Wiles
(Theorem A, Section 4.5 in~\cite{SW-reducible}).
\begin{theorem}
If $\overline\rho$ is reducible, then $\rho$ is modular.
\end{theorem}
The only point to be noted is that if $\overline\rho$ is
reducible, then $\rho$ is necessarily ordinary at $\ell$, and so
the quoted result applies. This follows immediately from
Proposition 9.1.2 of \cite{B}, which shows that if
$\rho|_{D_\ell}$ is irreducible and crystalline then its reduction
is irreducible.  (See also
 \cite{D} for a similar argument in  the case $\ell =3$ , $w=1$ where $w$ is the non-zero
 Hodge-Tate weight.)

\begin{theorem}
Suppose that $\rho$ is ordinary at $\ell$, and also suppose that
$\overline\rho$ is irreducible and modular. Then $\rho$ is
modular.
\end{theorem}
This follows immediately from Theorem 5.2 of Skinner and Wiles
in~\cite{SW-ord}: since $\rho$ is ordinary at $\ell$ and has
determinant $\epsilon_\ell^3$, we have
$$\overline\rho|_{I_\ell}\sim \begin{pmatrix}
\epsilon_\ell^3\pmod{\ell}&*\\ 0 &1\end{pmatrix}$$ where $I_\ell$ is
the inertia group at $\ell .$ In particular, it is
$D_\ell$-distinguished (i.e., the two diagonal characters in the
residual representation are distinct). As we have started out with
the condition that $\overline\rho$ is modular, the conclusion
above is the one given by Theorem 5.2 of \cite{SW-ord}.

The above two results deal with the case when the representation
$\rho$ is ordinary at $\ell.$ For the case when $\rho$ is
crystalline but not ordinary, we have the following result of
Taylor,~\cite{taylor}:
\begin{theorem}Suppose that $\rho$ is not ordinary at $\ell$. Let
$F$ be a totally real field of even degree such that
$\textrm{Gal}(F/\Q)$ is solvable, and such that $\ell$ splits
completely in $F.$ Assume the following:
\begin{itemize}
\item There is a Hilbert modular form of weight 4 and level 1 on
$GL_2(\mathbb{A}_F)$ which gives $\overline\rho$.
\item $\overline\rho$ restricted to the absolute Galois group of
$F\left(\sqrt{(-1)^{(\ell-1)/2}\ell}\right)$ is irreducible.
\end{itemize}
Then $\rho$ is modular.
\end{theorem}
\begin{Proof} Using Theorem 2.6 and Theorem 3.2 of \cite{taylor},
we see that $\rho|_{G_F}$ is `modular'. Solvable base change
results of Langlands then enable us to deduce the modularity of
$\rho$ from that of $\rho|_{G_F}$. \end{Proof}

As for the auxiliary totally real $F$ that appears in the above
theorem, we have the following `level lowering' result of Skinner
and Wiles,~\cite{SW-level}.
\begin{theorem} If $\overline\rho$ is irreducible and modular,
then we can find  a totally real solvable extension $F$, and a
Hilbert modular form of weight 4, level 1  on $GL_2(\mathbb{A}_F)$
which gives $\overline\rho$. (Note that $\overline\rho$ is
`crystalline'.) Moreover, given a finite set of primes $\Sigma$
where $\overline\rho$ is unramified, we can insist that the primes
in $\Sigma\cup\{\ell\}$ split completely in $F.$
\end{theorem}
\begin{Proof}
This is essentially proved in~\cite{SW-level} (they give a
complete proof for the case when the weight is 2).

Firstly, note that (cf.~\cite{dia-serre}) we can find a weight 4
newform $f$ on $\Gamma_0(N)$ with $N$ prime to $\ell$ and
divisible only by primes where $\overline\rho$ is ramified such
that
$$ \overline\rho\sim\rho_{f,\ell}\pmod{\ell}.$$
 We can therefore find an even degree,
totally real, solvable extension $F_1/\Q$ in which all the primes
in $\Sigma\cup\{\ell\}$ split completely, and a Hilbert modular
form $f_1$ over $F_1$ of level $\mathfrak{n}$, weight 4, such that
\begin{itemize}
\item $\overline\rho_{f_1,\ell}\sim \overline\rho|_{G_{F_1}},$
\item $\mathfrak{n}$ is an ideal dividing $N$ (and so is
prime to $\ell$), and
\item $\overline\rho|_{G_{F_1}}$ is unramified at all finite
places not dividing $\ell.$
\end{itemize}
Here, we have used the fact one can find solvable extensions with
prescribed local extensions (see for instance Lemma~2.2
of~\cite{taylorartin}).

It is now immediate that the result we want follows immediately
from a higher weight version of the main theorem of
\cite{SW-level}. The proof of the main theorem of \cite{SW-level}
for the higher weight case works in exactly the same way as for
the weight 2 case. The only modification, or something that is to
be checked, required is Lemma 2 of \cite{SW-level}; and, this can
be found in \cite{Ta1} (with a slight modification,
 the result we need is Lemma 2, {\it ibid.}).
\end{Proof}

\subsection{}

We now turn to the other serious issue of when the residual
representation is modular. Conjecturally, the answer is always yes
for odd, absolutely irreducible, continuous representations
(Serre's conjecture); our knowledge covers only a handful of
residue fields.

\begin{theorem} Let $\overline\rho :G_\Q\longrightarrow
GL_2(\F_5)$ be an absolutely irreducible,
 continuous representation with determinant equal
to $\epsilon_5^3\pmod{5}$. Then $\overline\rho$ is modular.
\end{theorem}
\begin{Proof} Let $\widetilde\rho:=\overline\rho \otimes
\epsilon_5^{-1}\pmod{5}.$ Then $\widetilde\rho$ is a continuous
Galois representation with values in $GL_2(\F_5)$ and determinant
equal to the mod 5 cyclotomic character $\epsilon_5\pmod{5}$. By
Theorem 1.2 of \cite{sbt}, we can find an elliptic curve $E$ over
$\Q$ such that $\overline\rho_{E,5}$, the Galois representation
given by the 5-torsion points of $E$, is equivalent to
$\widetilde\rho.$ Since all elliptic curves over $\Q$ are modular
(by \cite{bcdt}), it follows that $\widetilde\rho$ is modular. The
required result then follows.\end{Proof}

For Galois representations with values in $GL_2(\F_7)$, we have
the following (see \cite{mano}, and Theorem 9.1 of \cite{jama}):
\begin{theorem}
Let $\overline\rho:G_\Q\longrightarrow GL_2(\mathbb{F}_7)$ be an
absolutely irreducible, continuous, odd representation. If the
projective image of $\overline\rho$ is insoluble, we also assume
that:
\begin{itemize}
\item The projective image of inertia at $3$ has odd order.
\item The determinant of $\overline\rho$ restricted to the inertia
group at $7$ has even order.
\end{itemize}
Then $\overline\rho$ is modular.
\end{theorem}
The result applies, for example, if $\overline\rho$ has
determinant $\epsilon_7^3\pmod{7}$ and is unramified at 3.

The two results above do not specify the weight, or the level. For
the applications we have in mind, one needs to be able to produce
newforms of the `right weight and right level' giving rise to the
residual representation $\overline\rho.$ In view of
\cite{SW-reducible}, we might as well assume that $\overline\rho$
is absolutely irreducible. And once we make this assumption, we
can indeed produce newforms of the `right weight and right level'
(see Theorem 1.1 of~\cite{dia-serre}).

\section{Proof of the main theorem}

We now apply the results of the last section to the  given rigid
Calabi-Yau threefold $X$ defined over $\Q$.

\subsection*{X has good reduction at 3 and 7}
We consider the $7$-adic representation
$$\rho_{X,7}:G_\Q\longrightarrow GL_2(\Q_7).$$
Note that it is crystalline at $7$, and unramified at $3$. We also
know that it is absolutely irreducible.

If the residual representation  is absolutely  reducible, then
$\rho_{X,7}$ is modular by Theorem 2.2. So suppose that the
residual representation is absolutely irreducible. The residual
 representation is modular
by Theorem 2.7. Moreover, we know (from \cite{dia-serre}) that
there is weight $4$ newform of level coprime to $7$ and $3$ which
gives rise to $\overline\rho.$ The theorem in the ordinary case
follows from Theorem~2.3.

Finally, suppose that we are in the crystalline but not ordinary
case. Let $\omega_2$ denote the second fundamental character of
the inertia group at $7$. We then have (cf. \cite{B}, Chapter 9)
$$\overline\rho_{X,7}|_{I_7}\sim \left(
\omega_2^3\oplus\omega_2^{7.3}\right),$$ and so the projectivization of the
  image of
inertia is a cyclic group of order 8.  One can then check the
following:
  $\overline\rho_{X,7}$ restricted to $\mathbb{Q} \left(\sqrt{-7} \right)$ is
irreducible. Making use of the  (essential) fact that we are
working with a prime as large as $7$, the description above of
$\overline\rho_{X,7}|_{I_7}$
  and the inequality $(\ell + 1) / \gcd(w , \ell + 1) > 2 $ valid for $\ell = 7$ and
  $w=3$ ($w$ denotes the non-zero Hodge-Tate weight), one deduces the absolute
  irreducibility
   of $\overline\rho$ when restricted to the Galois group of
$\Q_7(\sqrt{-7})$. The proof (for general $\ell$ and $w$) follows
an idea of Ribet to deal with dihedral primes for
 Galois representations, and is given in full detail in \cite{D}
for the case $\ell = 3, \; w = 1$. Alternatively, one can
manipulate in the group $PSL_2(\F_7)$ (which we leave for the
interested reader).

Now take $F$ to be a totally real, even degree, solvable extension
of $\Q$ in which $7$ splits completely and $\overline\rho|_{G_F}$
arises from a Hilbert modular form of level 1. This can be done by
Theorem 2.5. Finally, we deduce that $\rho$ is modular by applying
Theorem 2.4 (that $\overline\rho$ restricted to the absolute
Galois group of $F(\sqrt{-7})$ is irreducible follows from the
preceding paragraph as $7$ splits completely in $F$).

\subsection*{X has good reduction at 5}

In this case, we make use of the $5$-adic representation
$\rho_{X,5}$. As before, we know that $\rho_{X,5}$ is absolutely
irreducible. Moreover, as before, we only need to consider the
cases when $\overline\rho_{X,5}$ is absolutely irreducible. Using
Theorem 2.6 and Theorem 1.1 of~\cite{dia-serre}, we see that there is
a newform of weight $4$ and level coprime to $5$ which gives rise
to $\overline\rho.$ One then deals with the ordinary case as
before using Theorem 2.3.

So suppose that we are in the crystalline but not ordinary case.
The image of inertia  at $5$ is then cyclic of order $8$; the
image of a decomposition group at $5$ is non-abelian of order
$16.$ Thus if $F$ is a totally real field in which $5$ splits
completely, we will still have
$$\left.\overline\rho_{X,5}\right|_{G_F}$$ absolutely irreducible.

Suppose that the restriction of~$\overline\rho_{X,5}$ to the
Galois group of~$G_{F(\sqrt{5})}$ is absolutely reducible.
Since~$F(\sqrt{5})$ is totally real, this is the same as
reducibility. In particular, we can assume that the image
of~$\overline\rho_{X,5}(G_{F(\sqrt{5})})$ is the subgroup
of~$GL_2(\F_5)$ consisting of elements
$$\begin{pmatrix} x &0\\ 0 & \pm x^{-1}\end{pmatrix}\quad
\textrm{with}\quad x\in\F_5^\times,$$ and hence that the image
of~$G_F$ is generated by the above subgroup and~$\begin{pmatrix}
0&1\\3&0\end{pmatrix}.$ It follows that elements with determinant
equal to~$2$ or~$3$ have trace~$0$.

We are now ready to make use of the condition that~$t_3(p)$ is not
divisible by~5. Choose an~$F$ which is totally real and solvable
as given by Theorem~2.5. We can choose such an~$F$ with~5 and~$p$
splitting completely in~$F$. Since~$t_3(p)$ is not divisible by~5,
the preceding argument shows  that~$F$ satisfies both the
conditions of Theorem~2.4. It follows that~$\rho_{X,5}$ is
modular.

\begin{remark}
We are grateful to Fred Diamond for pointing out that one could
prove modularity without appealing to base change. The relevant
result, dealing with modular representations over $\Q$, is proved
in \cite{DFG}.
\end{remark}

\begin{remark} What we proved, purely in terms of Galois
representations, is the following: \\
Let $\ell=5$ or $7,$ and let
$$\rho:G_\Q\longrightarrow GL_2(\Q_\ell^{ac})$$ be an irreducible continuous
representation, unramified outside finitely many primes, and
having determinant $\epsilon_\ell^3$. Assume the following:
\begin{enumerate}
\item  $\rho$ is crystalline at $\ell$ with Hodge-Tate
weights $\{0,3\}$.
\item The residual representation
$$\overline\rho:G_\Q\longrightarrow GL_2(\F_\ell^{ac})$$
takes values in $GL_2(\F_\ell).$
\item If $\ell=7$, then $\rho$ is unramified at $3$.
\item If $\ell=5$, then there is some prime $p\equiv \pm 2
\pmod{5}$ where $\rho$ is unramified and the trace of Frobenius at
$p$ is not divisible by 5.
\end{enumerate}
Then $\rho$ is modular.
\end{remark}

\section{Applications to known examples}

As mentioned in the introduction, all known rigid Calabi-Yau
threefolds are  modular. There are about thirty examples known;
the levels for which one can construct a rigid Calabi-Yau (of that
level) are
$$ 6,\ 8,\ 9,\ 12,\ 21,\ 25,\ 50,$$
 a family coming from conifolds (see Section~5.4 of~\cite{noriko-update}),
  and those of type $III_0$ (see Section~5.5 of~\cite{noriko-update})
 We are grateful to Helena Verrill and Noriko Yui for pointing these out; for
details and equations for the rigid Calabi-Yau threefolds,
see~\cite{noriko},~\cite{noriko-update} and the references there.
It is a pleasant surprise that the result of this article covers
all known rigid Calabi-Yau threefolds except for the level 9 case.

Let us check that the hypotheses of Theorem~1.2 are satisfied by
these known Calabi-Yau threefolds. In principle, we should write
down explicit models and verify the conditions  (for such explicit
models, see~\cite{noriko},~\cite{noriko-update}). However, we
already know that these examples are modular---we just want to
check that the result applies often. We get this information by
looking at tables of newforms. What we mean is the following: let
$N$ be a positive integer, and let $CY(N)$ be the set of rigid
Calabi-Yau threefolds whose $L$-series (or $q$-series) expansions
agree with newforms of weight~4, level $N,$ at least up to, say,
twenty terms. We know  from explicit calculations that if $X\in
CY(N),$ then $X$ has good reduction outside primes dividing $N.$
To see if our result would apply to the threefolds in $CY(N)$, we
need to see if either $(21,N)=1$ or $5\not|N.$ In the second case,
we would need to look at a few $t_3(p)$: we know how to calculate
them, and they agree with what comes from modular forms.

 For example, Theorem 1.2 with the first
condition will apply to the cases for which the level comes out to
be 8, 25 or 50. This is easy to check: for instance, the quintic
of Chad Schoen given in section 5.1 of \cite{noriko} is easily
seen to have good reduction at 3 and 7.
 For the cases when the level is 6,
12, or 21, one looks at tables of newforms with weight 4 (see
W.~Stein's website,~\cite{stein}). We need only look for newforms
with $q$-expansions having $\mathbb{Z}$-coefficients. We find
 one newform of level 6 (with $t_3(7)=-16$), one of level 12 (with
 $t_3(7)=8$), and two of level 21 (they have $t_3(13)=-34, \ -62$
 respectively); in all cases, we can apply Theorem 1.2 with the
 second set of conditions (the true values of these traces for the
 corresponding examples of rigid Calabi-Yau threefolds have been
 computed by many people and of course they agree with the above
 values).

We leave it to the reader to apply the result and check for
modularity of the family coming from conifolds and those of type
$III_0$ using the information given in Section~5.4
of~\cite{noriko-update}.

\begin{remark} Of course, our theorem does not provide any new
examples of rigid Calabi-Yau threefolds over $\Q.$ However, if
more examples are found, our theorem seems simpler to apply than
other results, and numerical evidence suggests that the hypotheses
seem to be satisfied sufficiently often for the theorem to be
genuinely useful. If an infinite family of rigid Calabi-Yau
threefolds over $\Q$ were to be discovered, the theorem might
apply to establish modularity for the family.
\end{remark}

\noindent{\sc Current addresses:}

\medskip

\noindent
L.~Dieulefait,\\
 \'{E}quipe de Th\'{e}orie des Nombres,
 Institut de Math\'{e}matiques de Jussieu,
Universit\'{e} Pierre et Marie Curie (Paris VI), Paris, France.\\
e-mail: {\tt luisd@math.jussieu.fr; ldieulefait@crm.es}
\medskip

\noindent
J.~Manoharmayum,\\
 Department of Pure Mathematics, University of
Sheffield, Sheffield S3 7RH, U.K.\\
e-mail: {\tt j.manoharmayum@sheffield.ac.uk}


\begin{thebibliography}{99}

\bibitem{Ber} P. Berthelot, {\it Alt\'{e}rations de vari\'{e}t\'{e}s alg\'{e}briques
[d'apr\`{e}s A. J. de
 Jong], } in ``S\'{e}minaire Bourbaki", vol. 1995/96, exp. 815, 273--311; Ast\'{e}risque
 1997.

\bibitem{B} C. Breuil, {\it $p$-adic Hodge theory, deformations and local
Langlands,}
notes of a course held at the CRM, Bellaterra, Spain, July 2001.


\bibitem{bcdt} C.~Breuil, B.~Conrad, F.~Diamond and
R.~Taylor, {\it On the modularity of elliptic curves over
$\mathbb{Q}$: wild 3-adic exercises, }  J. Amer. Math. Soc. {\bf
14} (2001), no. 4, pp 843--939.


\bibitem{dia-serre} F.~Diamond, {\it The refined conjecture of
Serre,} Elliptic curves, modular forms, $\&$ Fermat's last theorem
(Hong Kong, 1993), 22--37, Ser. Number Theory, I, Internat. Press,
Cambridge, MA, 1995.

\bibitem{DFG} F.~Diamond, M.~Flach and L.~Guo, {\it Adjoint motives of modular
 forms and the Tamagawa number conjecture,} preprint.


\bibitem{D} L. Dieulefait, {\it Modularity of abelian surfaces with quaternionic
 multiplication,} Math. Res. Lett., {\bf 10}
(2003), no. 2-3.

 \bibitem{FM} J. M. Fontaine and W. Messing, {\it $p$-adic periods and $p$-adic
\'{e}tale cohomology}, in ``Currents Trends in Arithmetical
Algebraic Geometry (Arcata, Calif., 1985)". Contemporary
Mathematics {\bf 67} (1987) 179--207, American Mathematical
Society.


\bibitem{jama} A.~F.~Jarvis and J.~Manoharmayum, {\it Modularity
of elliptic curves over totally real fields,} preprint.


\bibitem{mano} J.~Manoharmayum, {\it Modularity of certain
$GL_2(\mathbb{F}_7)$ representations,} Math. Res. Lett., {\bf 8}
(2001), no. 5-6, 703--712.

\bibitem{sbt} N.~I.~Shepherd-Barron, R.~Taylor, {\it Mod 2
and mod 5 icosahedral representations}, {\rm J. Amer. Math. Soc.}
{\bf 10} (1997) 283--298


\bibitem{SW-reducible} C.Skinner, A.Wiles, {\it Residually reducible representations
and modular forms}, Publ. Math. IHES 89 (2000) 5--126.

\bibitem{SW-ord} C.~Skinner, A.~Wiles, {\it Nearly ordinary deformations of
irreducible residual representations}, Ann. Fac. Sci. Toulouse
Math. (6) {\bf 10} (2001) 185--215


\bibitem{SW-level} C.Skinner, A.Wiles, {\it Base change and a problem of Serre},
Duke Math. J. 107 (2001) 15--25.

\bibitem{stein} W.~Stein,  {\tt
http://modular.fas.harvard.edu/Tables/index.html}

\bibitem{Ta1} R.Taylor, {\it On Galois representations associated to Hilbert
modular forms}, Invent. Math. 98 (1989) 265--280.

\bibitem{taylorartin} R.~Taylor, {\it On icosahedral Artin representations II},
to appear in Amer. J. Math. Available at {\tt
http://www.math.harvard.edu/\~{}rtaylor}

\bibitem{taylor} R.Taylor, {\it On the meromorphic continuation of degree two
L-functions}, preprint. Available at {\tt
http://www.math.harvard.edu/\~{}rtaylor}.

\bibitem{noriko} N.~Yui, {\it Arithmetic of certain Calabi-Yau varieties
and mirror symmetry,} Arithmetic algebraic geometry (Park City,
UT, 1999), 507--569, IAS/Park City Math. Ser., 9, Amer. Math.
Soc., Providence, RI, 2001.



\bibitem{noriko-update} N.~Yui, {\it Update on the Modularity of
Calabi-Yau Varieties,} this volume.


\end{thebibliography}
\end{document}